# Wordle: A Microcosm of Life. Luck, Skill, Cheating, Loyalty, and Influence!

*James P. Dilger*

**Abstract**

Wordle is a popular, online word game offered by the New York Times (nytimes.com). Currently there are some 2 million players of the English version worldwide. Players have 6 attempts to guess the daily word (target word) and after each attempt, the player receives color-coded information about the correctness and position of each letter in the guess. After either a successful completion of the puzzle or the final unsuccessful attempt, software can assess the player's luck and skill using Information Theory and can display data for the first, second, ..., sixth guesses of a random sample of all players. Recently, I discovered that the latter data is presented in a format that can easily be copied and pasted into a spreadsheet. I compiled data on Wordle players' first guesses from May 2023 - August 2023 and inferred some interesting information about Wordle players. A) Every day, about 0.2-0.5% of players solve the puzzle in one attempt. Because the odds of guessing the one of 2,315 possible target words at random is 0.043%, this implies that 4,000 - 10,000 players cheat by obtaining the target word outside of playing the game! B) At least 1/3 of the players have a favorite starting word, or cycle through several. And even though players should be aware that target words are never repeated, most players appear to remain loyal to their starting word even after its appearance as a target word. C) On August 15, 2023, about 30,000 players abruptly changed their starting word, presumably based on a crossword puzzle clue! Wordle players can be influenced! This study goes beyond social media postings, surveys, and Google Trends to provide solid, quantitative evidence about cheating in Wordle.

**Introduction**

Wordle began in June 2021 as a private game among the developer, Josh Wardle, and his family and friends. Five months later, Wardle made the game public and in January 2022, it was purchased by the New York Times. There are now millions of players of the English language version and there are versions in at least 50 other languages(*Wikipedia Entry: Wordle*, 2023)! The New York Times version is free to all players.

Every day, a new 5-letter target word is chosen (without replacement) from among 2,309 possible words. Players must guess the target word within 6 attempts, by guessing one of 12,545 "acceptable" English 5-letter words (Broz, 2023). After each attempt, the player receives color-coded information about the correctness and position of each letter in the guess. In the "normal" version of Wordle, players may use any "acceptable" word as their subsequent guess. However, in the "hard" mode, players are required to use information obtained in prior guesses when they make their subsequent guesses. At the conclusion of the game, the software reports the player's score as the number of guesses required and reports the player's statistics for the number of games played.

If the player is a New York Times subscriber, they may use the Wordlebot (Katz & Conlon, 2022a, 2022b). This software computes "skill" and "luck" scores for the player and compares these to the average of all players that day. The software uses Information Theory for these calculations (3Blue1Brown, 2022)

WordleBot also reports some data about the Share Percent and Rank of first, second, … sixth guessed words. The purpose of this paper is to try to infer some behaviors of Wordle players based on the Share Percent and Rank of first guessed words.

**Methods**

For May 3 – August 31, 2023, I used the "Compare and view your recent scores" option of WordleBot 2.0 (Katz & Conlon, 2022b) (Fig 1A) to view the analyses of past Wordles. According to WordleBot, the data represent a random sampling of all the approximately 1.7 million submissions that day. After the individual analysis, there is a screen "Everyone's Guesses Through Turn 1" (Fig 1B). This is a bar graph showing the top 15 first guesses of players, and 5 additional guesses that include your own first guess (TEACH, in my case), the target word itself (QUEST, in this case) and a few other random words. I thought it would be interesting to have this information in spreadsheet form, but I wasn't about to copy it manually. On a whim, I tried "Edit > Select All" and "Edit > Copy" and pasted the clipboard contents into Excel. The surprising result is shown in Fig 1C; ascii data in a single column format! After applying a few simple Excel formulas and operations, I had the entire list in 3 columns, sorted by rank (Fig 1D). And as if one surprise a day was not enough, I saw that a total of 50 first guesses were ranked! There is no indication on the Wordle website that any of this is possible.

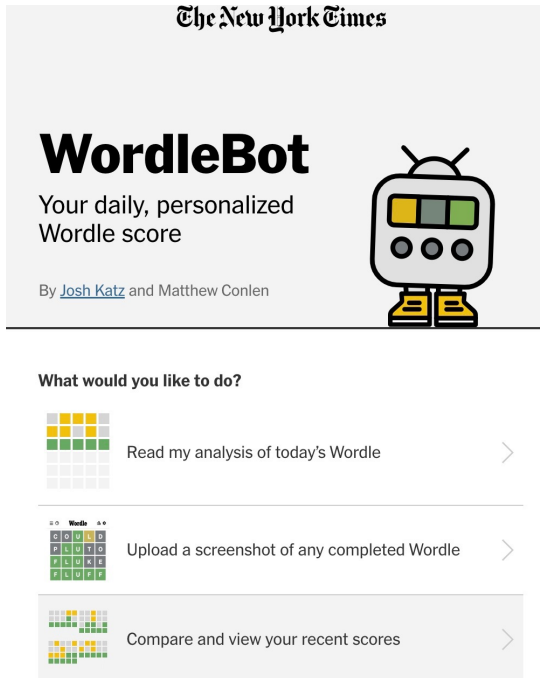
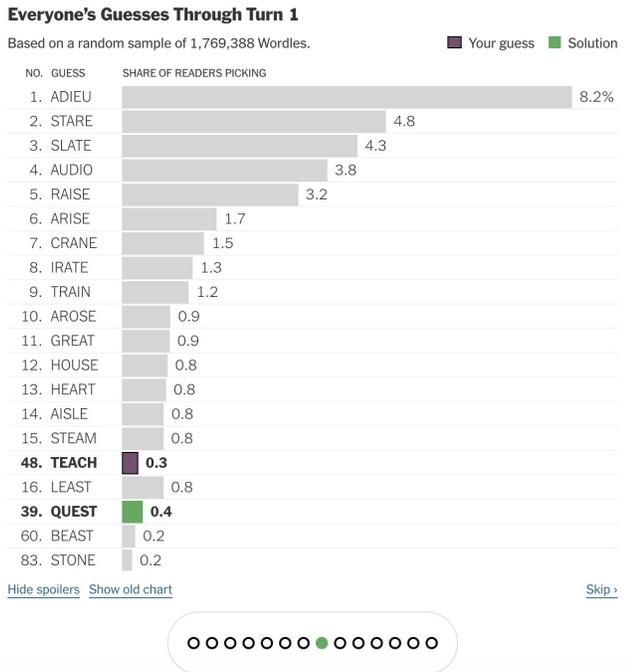

*Figure 1. A-D, the steps used to extract starting word information from WordleBot. "Select All" and "Copy" is used at step B. "Paste" into a spreadsheet is used at step C. Results are put into columns and sorted in step D. The target word of the day, QUEST, is in red font.*

I compiled the first guess data in a Microsoft® Excel for Mac spreadsheet (v.16.75.2). I used Igor Pro 8 (v.8.04, WaveMetrics, Inc., Lake Oswego, Ore) for graphing and analysis. Two measures of a starting word's popularity, Rank and Share Percent of readers are provided. Share Percent is quantized in 0.1% increments and low values are denoted as "<0.1"; it is a useful metric for first guess words with high popularity. Conversely, Rank works better with less popular first guess words. Because of possible confusion with the concept of Rank (where smaller numbers are more popular), I include a "thumbs-up" or "thumbs-down" icon at the top of all ordinate axes to indicate "better / more popular" and "worse / less popular" in all figures.

The data come from Wordle played in the "normal" mode, currently about 1.7 million players do this. In the "hard" mode, players must use all revealed hints in subsequent guesses. About 220,000 people play in this mode. A few results from "hard" Wordle are presented at the end of the Results section.

**Results**

The probability of guessing the target word by chance is 1/2,315 = 0.043% (assuming that the player "knows" the list of target words). That's 860 players. But look at Figure 2A which shows 4 months of first guess data in terms of Share Percent. Do I mean to tell you that never, not once, was the Share Percent of the first guess less than 0.2% (4,000 players)? Yup! So, at least 5-times as many people as expected were "guessing" the target word by "chance"! Okay, maybe we should consider that some 800 target words have already been used (and will not be re-used), the lucky guesser would still have only 1/1,515 = 0.066% odds. Yet, it happens consistently every day! Some days it's as high as 0.5% (10,000 players). What shall we call these people? Hmmm, "cheaters" comes to mind, so that's what I call 'em!

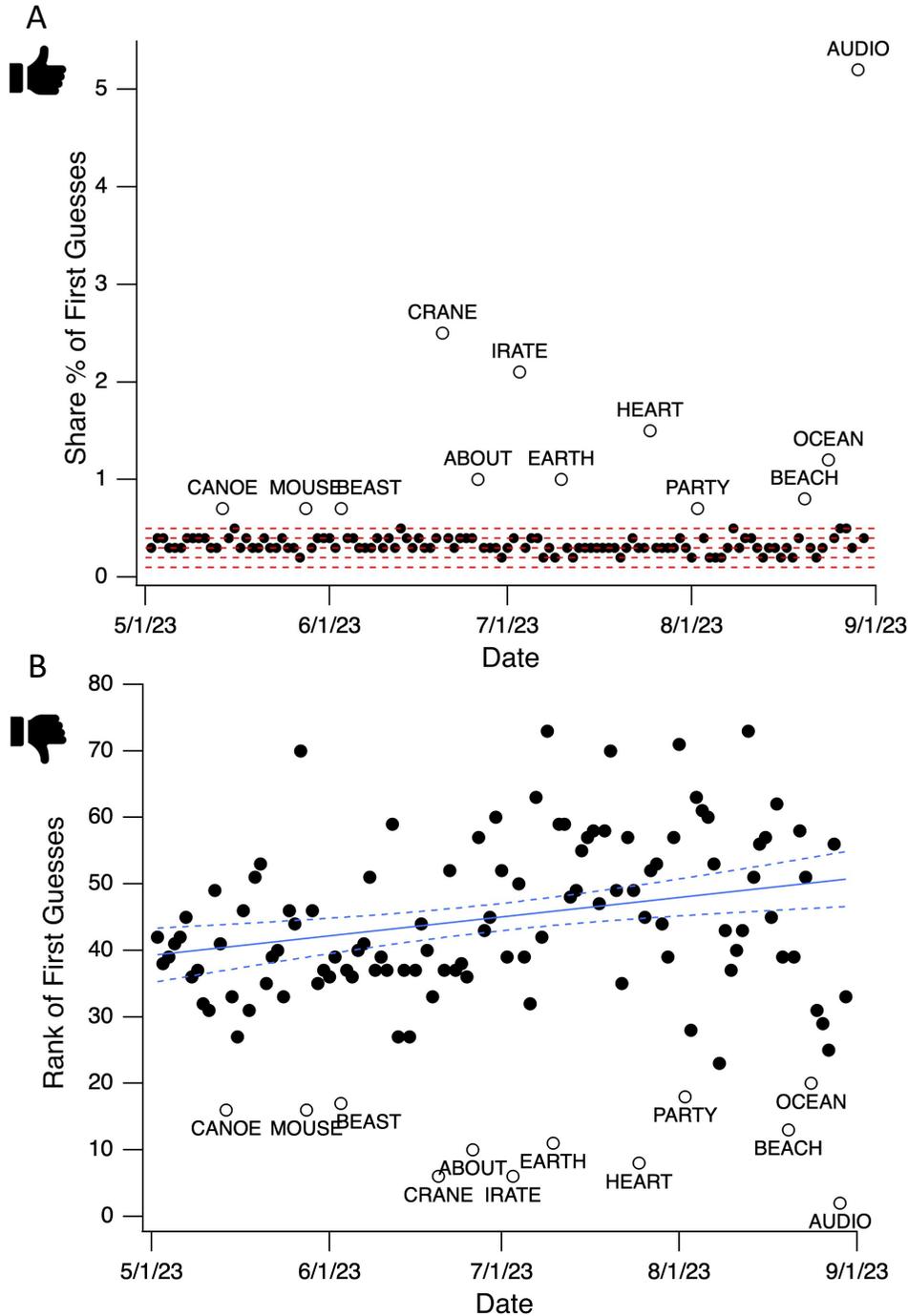

*Figure 2. A. Share Percent of all First Word Guesses. Every day, many players "guess" the target word on the first attempt; many more than expected from chance (0.04%). Low ranking numbers correspond to popular guesses. The target word ranks between the 20th and 70th most popular first guess. The labeled open symbols correspond to words that are popular first guesses on most days. The blue lines are the best fit linear regression line and 95% confidence limits.*

What should we make of the days when 1.0 – 5.2% of players guess right the first time?

Well, it could be that the target word on those days was a popular first guess on almost any day.

This is seen in Figure 2B where the data are expressed as Rank and such words are shown as open circles and labeled. CANOE, for example, had a rank of 60 (0.2% share) on a random day and this became a rank of 16 (0.7% share) on its special day. EARTH is a very popular starting word; it consistently ranks in the 20s and gets a 0.5% share. So, it isn't too surprising that it ranked 11 (1% share) on EARTH day! For both CANOE and EARTH, the jump in share probably represents our old friends, the cheaters!

Omitting the "normally popular" words (Rank≤20), I performed linear regression on Rank vs Date and found an upward trend (Fig 2b). Does this mean that cheating is on the wane? It certainly looks that way; the slope of the line is significantly different from zero: 0.095 ± 0.06 rank level per day (95% confidence limits). Are cheaters getting bored with Wordle, or have they just been on summer vacation?

Just to hammer down the issue with cheating, I looked for cases where there was ranking data for target words, prior to their being chosen as target words (omitting the normally popular first guesses). This wasn't always possible, but the results are compelling (Fig. 3). Whereas the Rank of every target word was <80 on its special day, it was, at best, >200 on a previous day. What serious Wordle player would choose NANNY as a first guess? You'd be testing only 2 vowels and 1 consonant! And IGLOO? srsly?

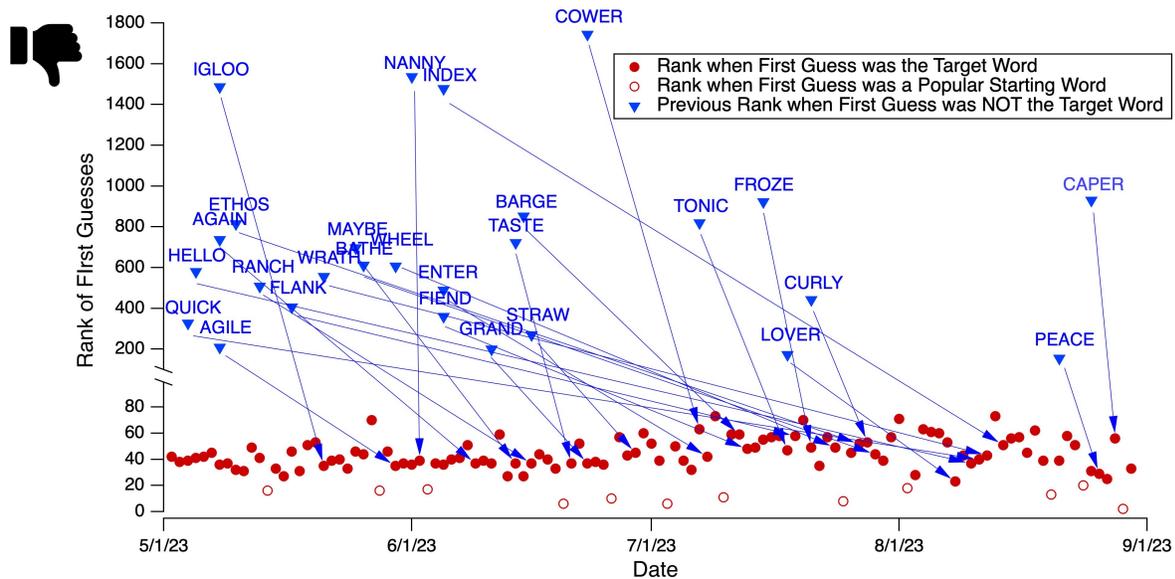

Figure 3. Rank of First Guesses. The lower portion of the graph is a repeat of Figure 2B. The upper portion shows instances when a word was guessed on a day prior to its being the target word.

Enough about the cheaters! Well, almost. Let's talk about starting word loyalty. My personal starting word is TEACH. I realize that it isn't the strongest starting word, but it's my word and I'm sticking with it! I'm loyal! It seems that there may be a cadre of loyal TEACH fans out there because TEACH consistently ranks in the 40s or 50s with a 0.3% share. (Fig 4). In fact, it appears to be gaining in popularity. Linear regression of Rank vs Time shows a significant trend towards lower ranking numbers (more popular) (-0.07±0.03 rank level per day, p<0.001).

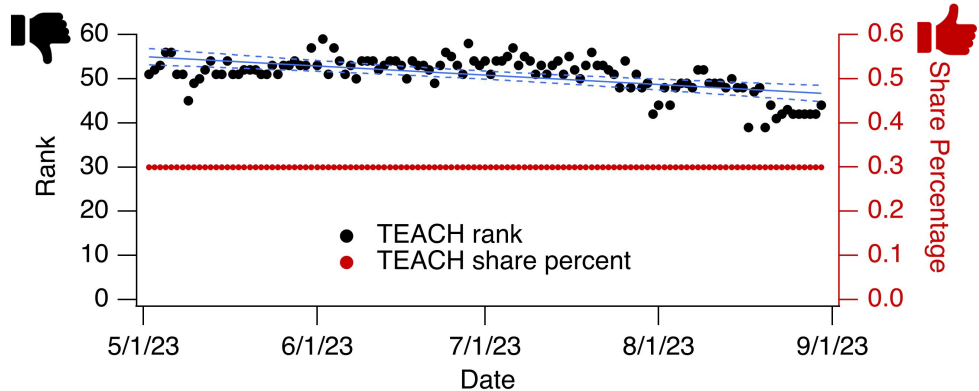

*Figure 4. Rank and Share Percent for TEACH as a starting guess. TEACH seems to have a loyal following. its popularity is growing as seen by the significantly nonzero slope of the linear regression line. Dashed lines indicate 99.9% confidence intervals.*

It is interesting to look at some of the very popular first guesses (the top 15 are always provided by WordleBot) especially when one happens to be the target word of the day. Fig 5 shows the four cases where this occurred during the data collection time: AUDIO typically ranks 4, CRANE typically ranks 6 or 7, IRATE typically ranks 8 or 9, and HEART typically ranks 10 or 11. All four of them garnered a steady Share Percent until … the day when the word was the target word! On that day, the Share Percent spiked (due to the "You-Know-Whos") and thereafter the Share Percent fell by about 10-15% of the previous values and stabilized there (you'll have to be patient to see how this pans out for AUDIO). It looks as though most of the loyal first guess players remain loyal even after the word has been the target! After TEACH is the target word, I plan to stay loyal as well! (Maybe I am just too lazy to rethink my strategy). Another measure of loyalty is the fact that the top 15 first words consistently have a combined Share Percent of about 30% (data not shown). Of course, this does not necessarily mean that these players always use the same first guess, they might rotate among several standard starting words.

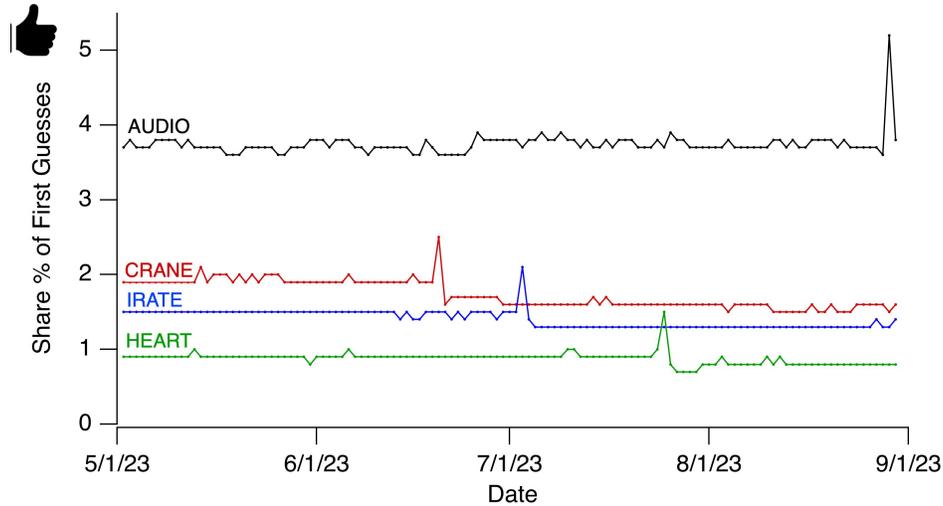

*Figure 5. Share Percent for four words that are popular first guesses and happened to be a target word during the data collection period. Each word exhibits a spike in popularity on the day of its being the target word (thanks to cheaters). Thereafter, some players abandon the word as their first guess, but most remain loyal.*

The number 1 starting "normal" mode word is consistently ADIEU with a Share of just under 8% (Fig 6). I suppose that the appeal of ADIEU is the fact that it assays 4 vowels in one shot. ADIEU has not been the target word yet, but it had a sudden surge in popularity on August 15. And so far, much of this change has persisted. What??? A likely explanation comes from observation that the August 15 NYT Mini Crossword clue for 6 across was "most popular starting guess in Wordle" and, of course, the answer was ADIEU. This seems to have inspired about 30,000 players to abandon whatever opening strategy they had been using and follow the crowd. ADIEU also spiked on August 30. Again, it was not the target word, but that was the day for AUDIO. Were players especially fond of words with 4 vowels that day?

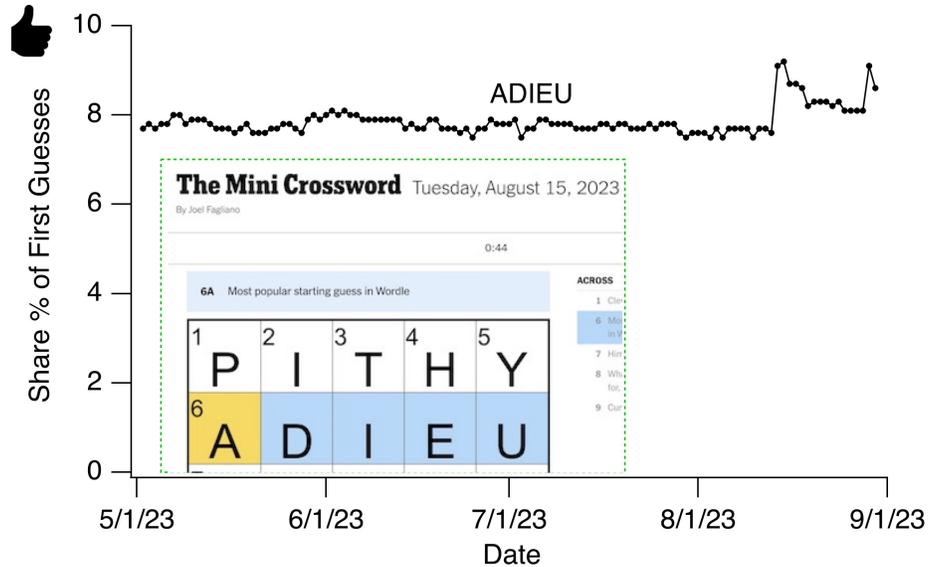

*Figure 6. The curious case of ADIEU. It has always been the top first guess. It became suddenly more popular on August 15, but not because ADIEU was the target word that day. Apparently, enthusiasts who play both Wordle and the Mini Crossword were inspired by a crossword puzzle clue! The spike on August 30 is more difficult to explain.*

With just one week left for data collection, I began to wonder about "hard" Wordle. I shelled out $4.07 for a new 4-week nytimes.com subscription (the things I do for Science!) and used a different browser (so I wouldn't have to constantly log in and out). Based on this small sample, "hard" core Wordle players cheat with their starting word a bit less than "normal" mode players (Table 1). The small differences might be attributed to roundoff errors in quantizing the Share Percent numbers, especially on days when the target word is not one of the most popular starting words. (If you are wondering why there is no data for 29-Aug, I'm really bummed that I forgot to do "hard" Wordle that day … so don't bring up the subject again!)

| Date | Word | Rank - normal | Rank - hard | Share - normal | Share - hard |
|---|---|---|---|---|---|
| 25-Aug | OCEAN | 10 | 12 | 1.2 | 0.8 |
| 26-Aug | CHOIR | 31 | 36 | 0.4 | 0.3 |
| 27-Aug | PEACE | 29 | 45 | 0.5 | 0.3 |
| 28-Aug | WRITE | 25 | 33 | 0.5 | 0.4 |
| 29-Aug | CAPER | 56 |  | 0.3 |  |
| 30-Aug | AUDIO | 2 | 5 | 5.2 | 4.0 |
| 31-Aug | BRIDE | 33 | 45 | 0.4 | 0.3 |

**Discussion**

Cheating is a fact of life. It's probably even built into our genes! (Sun, 2023) With multi-player video games, "to advance towards completion" and "gain advantage over others" were cited as being common justifications for cheating (Doherty et al, 2014). But why should someone cheat at a single-player, inconsequential game? A good question, indeed! Cheating in one-player games has been examined by others using different techniques. In a pre-Wordle conference proceedings (Passmore et al, 2020), a survey was used to assess 188 peoples' level of cheating (aka using "extraneous game advantages") and their feelings about cheating. Many players became "frustrated" at some point in the game and then felt "joy" or "relief" after having surpassed the hurdle with a cheat. Cheating at Wordle *per se* has been addressed by several groups. One group used Google Trends to identify searches for the Wordle target word (Wormley and Cohen, 2022). They were not able to quantify cheating this way. Instead, they calculated a "Wordle cheating index" for each US State and looked for correlations between the index and "religiosity" and "cultural tightness". Needless to say, this generated some state-pride issues, and I just don't want to go down that Reddit hole! Google Trends data was also used to propose that Wordle cheating rose steadily in early 2022 (Iancu, 2022).

Some quantitative data for Wordle cheating does exist. Only 0.02% of players who post their results on social media, report being successful on their first try (Broz, 2023). This is lower than expected from chance, so may indicate that some players start with a word that is not one of the possible target words. On the other hand, a survey of 1,000 players revealed that >10% admitted to cheating (solitaired.com, 2022). Of course, that number includes more than just those who cheat on the first try. I would surmise that players become more likely to cheat as they use up their allotted 6 attempts.

The New York Times provided some information about player's starting words (Amlen, 2022). Looking over the early summer 2022 timeframe, the five most popular starting words were ADIEU, AUDIO, STARE, RAISE, and ARISE. In my observations, these five were in the top 6 of most frequently used starting words one year later. Amlen estimated that 28% of players use the same, or a few, starting words consistently.

This study is the first one to go beyond social media postings, surveys, and Google Trends to provide solid, quantitative evidence about cheating in Wordle. It is naturally restricted to cheating on the first guess and is limited to the time period during which I was able to use the WordleBot to view historical data (currently restricted to 90 days).

**Conclusions**

Wordle is known as a game of luck and skill. What I set out to do here is to show that Wordle is much more than that. There is rampant cheating (just as in real life), loyalty (whether deserved or not), and influence (for good or evil, who knows?). It is also fun to play, but you must rely on my anecdotal evidence for that. I play daily, competing against my sister and daughter. It is a friendly competition. We cheer each other's successes and try to learn from our mistakes. Best of all, they tolerate my nerdy tendency for statistics and puns. We are baffled as to how first-word cheaters actually have fun playing but, that does not diminish our enjoyment of the game!

**Acknowledgements**

Thanks! To my Wordle partners Lena Corrigan and Emily Dilger for their resilience and encouragement. To my neighbor, workout partner and Operations Management expert Kats Sasanuma for his critical reading of the manuscript. To my best man Steve Crandall for publishing advice. To my loving wife Kate Dilger who often reminds me that there are uses for a big vocabulary other than for playing word games! (Uh, does being able to guess 5-letter words really constitute a big vocabulary?) To my son Andrew Dilger for his encouragement even though he cannot understand how I spend so much time with Wordle when I could be playing Chessle (https://jackli.gg/chessle/).